\documentclass{colt2016} 

\let\vec\undefined 
\usepackage{macrosArticle}
\usepackage{forest}
\usepackage[colorinlistoftodos, textwidth=26mm, shadow,color=yellow]{todonotes}

\def\UCB{\mathrm{U}}
\def\LCB{\mathrm{L}}
\def\Arms{\cP}

\def\ihat{\hat{\imath}}
\def\jhat{\hat{\jmath}}
\def\itilde{\tilde{\imath}}
\def\jtilde{\tilde{\jmath}}

\DeclareMathOperator{\KL}{KL}

\title[Maximin Action Identification]{Maximin Action Identification: A New Bandit Framework for Games}

 \coltauthor{\Name{Aur\'elien Garivier} \Email{aurelien.garivier@math.univ-toulouse.fr}\\
 	\addr Institut de Math\'ematiques de Toulouse; UMR5219\\
 	Universit\'e de Toulouse; CNRS\\
 	UPS IMT, F-31062 Toulouse Cedex 9, France
 	\AND
 	\Name{Emilie Kaufmann} \Email{emilie.kaufmann@inria.fr}\\
 	\addr Laboratoire CRIStaL ; Equipe SequeL \\
 	CNRS UMR ; Inria Lille - Nord Europe\\
 	59650 Villeneuve d'Ascq, France
 	\AND
 	\Name{Wouter M. Koolen} \Email{wmkoolen@cwi.nl}\\
 	\addr Centrum Wiskunde \& Informatica \\
 	Amsterdam, the Netherlands
 }

\begin{document}
\selectlanguage{english}

\maketitle

\begin{abstract}
We study an original problem of pure exploration in a strategic bandit model motivated by Monte Carlo Tree Search. It consists in identifying the best action in a game, when the player may sample random outcomes of sequentially chosen pairs of actions. 
We propose two strategies for the fixed-confidence setting: Maximin-LUCB, based on lower- and upper- confidence bounds; and  Maximin-Racing, which operates by successively eliminating the sub-optimal actions.
We discuss the sample complexity of both methods and compare their performance empirically. We sketch a lower bound analysis, and possible connections to an optimal algorithm. 

\end{abstract}

\begin{keywords}
multi-armed bandit problems, games, best-arm identification, racing, LUCB
\end{keywords}

\section{Setting: A Bandit Model for Two-Player Zero-Sum Random Games}

We study a statistical learning problem inspired by the design of computer opponents for playing games. We are thinking about two-player zero sum full information games like Checkers, Chess, Go \citep{deep.go} \dots, and also games with randomness and hidden information like Scrabble or Poker \citep{poker}. At each step during game play, the agent is presented with the current game configuration, and is tasked with figuring out which of the available moves to play. In most interesting games, an exhaustive search of the game tree is completely out of the question, even with smart pruning.

Given that we cannot consider all states, the question is where and how to spend our computational effort. A popular approach is based on Monte Carlo Tree Search (MCTS) \citep{go.review,SurveyMCTS12}. Very roughly, the idea of MCTS is to reason strategically about a tractable (say up to some depth) portion of the game tree rooted at the current configuration, and to use (randomized) heuristics to estimate values of states at the edge of the tractable area. One way to obtain such estimates is by `rollouts': playing reasonable random policies for both players against each other until the game ends and seeing who wins.

MCTS methods are currently applied very successfully in the construction of game playing agents and we are interested in understanding and characterizing the fundamental complexity of such approaches. The existing picture is still rather incomplete. For example, there is no precise characterization of the number of rollouts required to identify a close to optimal action. Sometimes, cumulated regret minimizing algorithms (e.g. UCB derivatives) are used, whereas only the simple regret is relevant here. As a first step in this direction, we investigate in this paper an idealized version of the MCTS problem for games, for which we develop a theory that leads to sample complexity guarantees.

\begin{figure}
\centering
\begin{tikzpicture}[every node/.style={minimum size=3em},node distance=2em and 1em]
\path
  node[circle,draw] (n) {$\max$}
  node[circle,draw,below left = 2em and 4em of n]    (n0) {$\min$}
  node[circle,draw, below right = 2em and 4em of n]  (n1) {$\min$}
  node[circle,draw,below left = of n0]   (n00) {$\mu_{1,1}$}
  node[circle,draw, below right = of n0] (n01) {$\mu_{1,2}$}
  node[circle,draw,below left = of n1]   (n10) {$\mu_{2,1}$}
  node[circle,draw, below right = of n1] (n11) {$\mu_{2,2}$};
\draw[->]
      (n) edge (n0)
      (n) edge (n1)
      (n0) edge (n00)
      (n0) edge (n01)
      (n1) edge (n10)
      (n1) edge (n11);
\end{tikzpicture}
\caption{Game tree when there are two actions by player $(K=K_1=K_2=2)$.}\label{fig:gametree}
\end{figure}
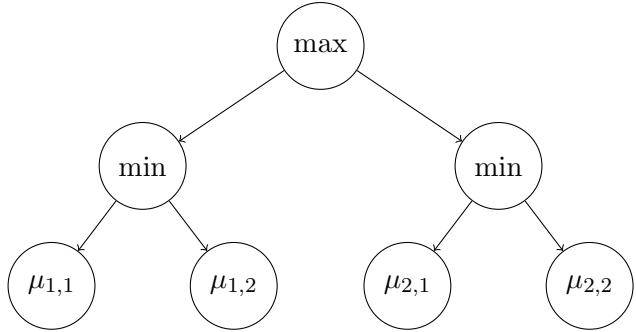

More precisely, we study perhaps the simplest model incorporating both strategic reasoning and exploration. We consider a two-player two-round zero-sum game, in which player A has $K$ available actions. For each of these actions, indexed by $i$, player B can then choose among $K_i$ possible actions, indexed by $j$.  For $i \in \{1,\dots,K\}$ and $j \in \{1,\dots,K_i\}$, when player A chooses action $i$ and then player B chooses action $j$, the probability that player A wins is $\mu_{i,j}$.
We investigate the situation (see Figure~\ref{fig:gametree} for an example) from the perspective of Player A, who wants to identify a maximin action
\[i^* \in \argmax{i \in \{1,\dots,K\}} \min_{j \in \{1,\dots,K_i\}} \ \mu_{i,j}.\]
Assuming that Player B is strategic and picks, whatever A's action $i$, the action $j$ minimizing $\mu_{i,j}$, this is the best choice for A.

The parameters of the game are unknown to player A, but he can repeatedly choose a pair $P=(i,j)$ of actions for him \emph{and player B}, and subsequently observe a sample from a Bernoulli distribution with mean $\mu_{i,j}$. At this point we imagine the sample could be generated e.g.\ by a single rollout estimate in an underlying longer game that we consider beyond tractable strategic consideration.
Note that, in this learning phase, Player A is not playing a game: he chooses actions for himself \emph{and} for his adversary, and observes the random outcome. 

The aim of this work is to propose a dynamic sampling strategy for Player A in order to minimize the total number of samples (i.e. rollouts) needed to identify $i^*$. Letting 
\[\cP = \big\{(i,j)  : 1 \leq i \leq K, 1 \leq j \leq K_i\big\},\]
we formulate the problem as the search of a particular arm in a stochastic bandit model with $\overline{K} = \sum_{i=1}^K K_i$ Bernoulli arms of respective expectations $\mu_{P}$, $P \in \cP$. In this bandit model, parametrized by $\bm \mu = (\mu_{P})_{P \in \cP}$, when the player chooses an arm (a pair of actions) $P_t$ at round $t$, he observes a sample $X_t$ drawn under a Bernoulli distribution with mean $\mu_{P_t}$. 

In contrast to best arm identification in bandit models (see, e.g., \cite{EvenDaral06,Bubeck10BestArm}), where the goal is to identify the arm(s) with highest mean, $\text{argmax}_P \ \mu_P$, here we want to identify as quickly as possible the maximin action $i^*$ defined above. 
For this purpose, we adopt a sequential learning strategy (or algorithm) $(P_t,\tau,\ihat)$. Denoting by $\cF_t = \sigma(X_1,\dots,X_t)$ the sigma-field generated by the observations made up to time $t$, this strategy is made of 
\begin{itemize}
 \item a sampling rule $P_t \in \cP$ indicating the arm chosen at round $t$, such that $P_t$ is $\cF_{t-1}$ measurable,
 \item a stopping rule $\tau$ after which a recommendation is to be made, which is a stopping time with respect to $\cF_t$,
 \item a final guess $\ihat$ for the maximin action $i^*$.
\end{itemize}
For some fixed $\epsilon \geq 0$, the goal is to find as quickly as possible an $\epsilon$-maximin action, with a high accuracy. More specifically, given $\delta \in ]0,1[$, the strategy should be \emph{$\delta$-PAC}, i.e.\ satisfy 
\begin{equation}\label{eq:PAC}
\forall \bm \mu, \ \bP_{\bm \mu} \left(\min_{j \in \{ 1 \dots K_{i^*}\}} \ \mu_{i^*,j} - \min_{j \in \{ 1 \dots K_{\ihat}\}} \ \mu_{\ihat,j} \leq \epsilon\right) \geq 1 - \delta,
\end{equation}
while keeping the total number of samples $\tau$ as small as possible. This is known, in the best-arm identification literature, as the \emph{fixed-confidence} setting; alternatively, one may  consider the \emph{fixed-budget} setting where the total number of samples $\tau$  is fixed in advance, and where the goal is to minimize the probability that $\ihat$ is not an $\epsilon$-maximin action.

\paragraph{Related work.}

Tools from the bandit literature have been used in MCTS for around a decade (see \cite{SurveyRemiMCTS} for a survey). Originally, MCTS was used to perform planning in Markov Decision Process (MDP), which is a slightly different setting with no adversary: when an action is chosen, the transition towards a new state and the reward observed are generated by some (unknown) random process. A popular approach, UCT \citep{KocsisBBMCP06} builds on Upper Confidence Bounds algorithms, that are useful tools for regret minimization in bandit models (e.g., \cite{Aueral02}). In this slightly different setup (see \cite{Bubeck:Survey12} for a survey), the goal is to maximize the sum of the sample collected during the interaction with the bandit, which amounts in our setting to favor rollouts for which player A won (which is not necessary in the learning phase). This situation is from a certain perspective a little puzzling and arguably confusing, because as shown by \cite{Bubeckal11}, regret minimization and best arm identification are incompatible objectives in the sense that no algorithm can simultaneously be optimal for both.

More recently, tools from the best-arm identification literature have been used by \cite{STOP14} in the context of planning in a Markov Decision Process with a generative model. The proposed algorithm builds on the UGapE algorithm of \cite{Gabillon:al12} to decide for which action new trajectories in the MDP starting from this action should be simulated. Just like a best arm identification algorithm is a building block for such more complex algorithms to perform planning in an MDP, we believe that understanding the maximin action identification problem is a key step towards more general algorithms in games, with provable sample complexity guarantees. For example, an algorithm for maximin action identification may be useful for planning in a competitive Markov Decision Processes \cite{Filar96CMDP} that models stochastic games. 

\paragraph{Contributions.} In this paper, we propose two algorithms for the maximin action identification in the fixed-confidence setting, inspired by the two dominant approaches used in best arm identification algorithms. The first algorithm, Maximin-LUCB, is described in Section~\ref{sec:MLUCB}: it relies on the use of Upper and Lower Confidence Bounds. The second, Maximin-Racing is described in Section~\ref{sec:MRacing}: it proceeds by successive eliminations of the sub-optimal arms. We prove that both algorithms are $\delta$-PAC, and give upper bounds on their sample complexity. Along the way, we also propose some perspectives of improvement that are illustrated empirically in Section~\ref{sec:Experiments}. Finally, we propose in  Section~\ref{sec:Perspectives} for the two-actions case a lower bound on the sample complexity of any $\delta$-PAC algorithm, and sketch a strategy that may be optimal with respect to this lower bound. Most proofs are deferred to the Appendix.

\paragraph{Notation.}

\begin{figure}[h]
\centering
\begin{tikzpicture}[every pin edge/.style={help lines,<-}]
\node (m11) at (0,0) [pin=below:$\mu_{1,1}$] {$\times$};
\node (m12) at (1,2) [pin=left:$\mu_{1,2}$] {$\times$};
\node (m13) at (2,3) [pin=left:$\mu_{1,3}$] {$\times$};
\draw[dashed] (3,-3.5) -- (3,4);
\node (m21) at (4,-1) [pin=below:$\mu_{2,1}$] {$\times$};
\node (m22) at (5,-.5)  [pin=right:$\mu_{2,2}$] {$\times$};
\node (m23) at (6,3.5)  [pin=left:$\mu_{2,3}$] {$\times$};
\draw[dashed] (7,-3.5) -- (7,4);
\node (m31) at (8,-3) [pin=right:$\mu_{3,1}$] {$\times$};
\node (m32) at (9,1)  [pin=left:$\mu_{3,2}$] {$\times$};
\node (m33) at (10,2)  [pin=left:$\mu_{3,3}$] {$\times$};
\draw [->] (m33) -- (m32);
\draw [->] (m32) -- (m31);
\draw [->] (m23) -- (m22);
\draw [->] (m22) -- (m21);
\draw [->] (m13) -- (m12);
\draw [->] (m12) -- (m11);
\draw [->] (m11) -- (m21);
\draw [->] (m21) -- (m31);
\end{tikzpicture}
\caption{\label{fig:ArmSetup}Example `normal form' mean configuration. Arrows point to smaller values.}
\end{figure}
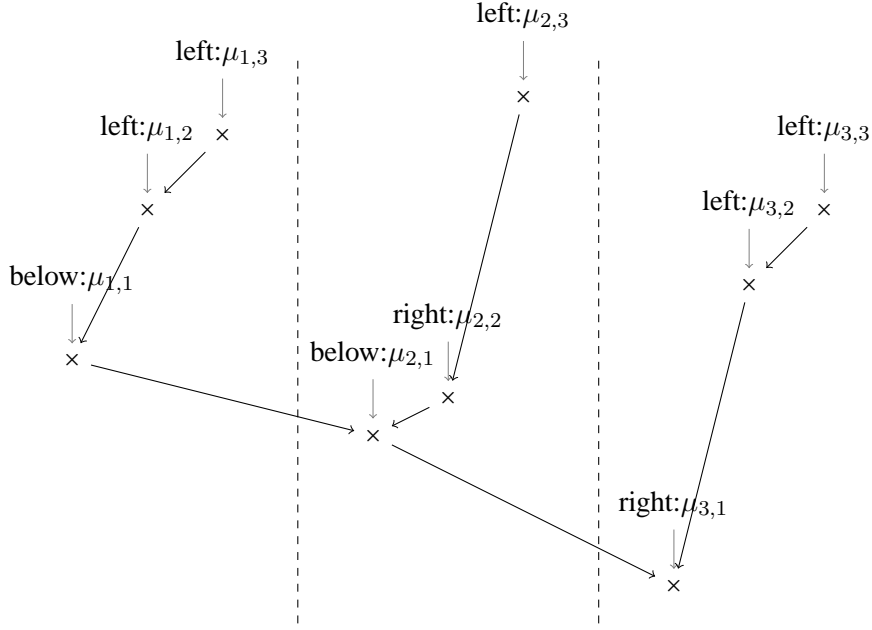

To ease the notation, in the rest of the paper 
we assume that the actions of the two players are re-ordered so that for each $i$, $\mu_{i,j}$ is increasing in $j$, and $\mu_{i,1}$ is decreasing in $i$ (so that $i^* = 1$ and $\mu^* = \mu_{1,1}$). These assumptions are illustrated in Figure~\ref{fig:ArmSetup}. With this notation, the action $\ihat$ is an $\epsilon$-maximin action if $\mu_{1,1} - \mu_{\ihat,1} \leq \epsilon$.
We also introduce $\Arms_i = \{(i,j), j \in \{1,\dots,K_i\}\}$ as the group of arms related to the choice of action $i$ for player A. 

\section{First Approach: M-LUCB}\label{sec:MLUCB}

We first describe a simple strategy based on confidence intervals, called Maximin-LUCB (M-LUCB). Confidence bounds have been successfully used for best-arm identification in the fixed-confidence setting (\cite{Shivaramal12,Gabillon:al12,Jamiesonal14LILUCB}). The algorithm proposed in this section for maximin action identification is inspired by the LUCB algorithm of \cite{Shivaramal12}, based on Upper and Lower Confidence Bounds.

For every pair of actions $P \in \cP$, let $\cI_P(t) = [\LCB_P(t),\UCB_P(t)]$ be a confidence interval on $\mu_P$ built using observations from arm $P$ gathered up to time $t$. Such a confidence interval can be obtained by using the number of draws $N_P(t) := \sum_{s=1}^t \ind_{(P_t = P)}$ and the empirical mean of the observations for this pair $\hat{\mu}_P(t) := \sum_{s=1}^t X_t \ind_{(P_{t}=P)}/N_P(t)$. 
The M-LUCB strategy aims at aligning the lower confidence bounds of arms that are in the same group $\Arms_i$. Arms to be drawn are chosen two by two: for any even time $t$, defining for every $i \in \{1,\dots,K\}$
\[c_{i}(t) = \argmin{1 \leq j \leq K_i} \ \LCB_{(i,j)}(t) \ \ \ \text{and} \ \ \ \hat{\imath}(t) = \argmax{i} \min_{j} \hat{\mu}_{i,j}(t)\;,\]
 the algorithm draws  at round $t+1$ and $t+2$ the arms
 \[H_t = (\hat{\imath}(t),c_{\hat{\imath}(t)}(t)) \ \ \ \text{and} \ \ \ S_t = \argmax{ P \in \{(i,c_i(t))\}_{i \neq \ihat}} \ \UCB_P(t).\]
This is indeed a regular LUCB sampling rule on a time-dependent set of arms each representing  one action: $\{(i,c_i(t))\}_{i \in \{1,\dots,K\}}$. In the two-actions case, one may alternatively draw at each time $t$ the arm $P_{t+1} = \text{argmax}_{P \in \{H_t,S_t\}} \ N_P(t)$ only. 
 
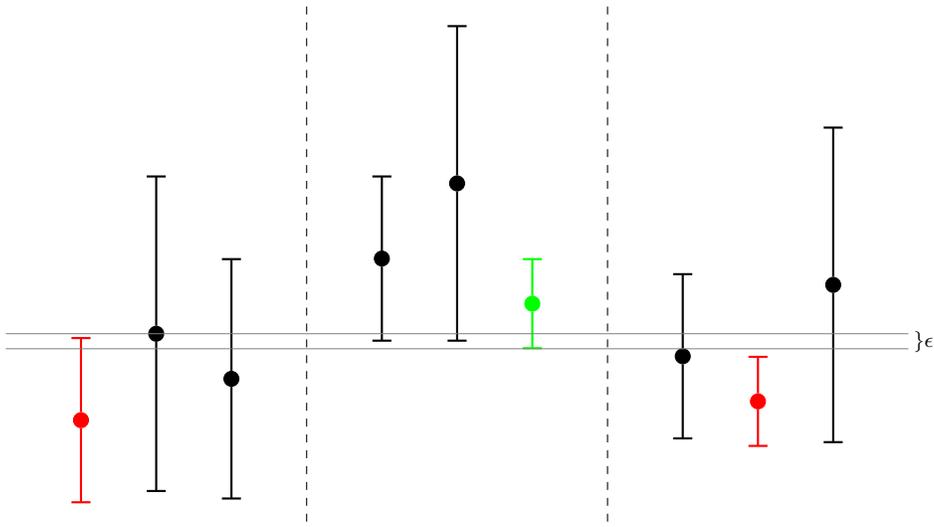
\begin{figure}
\centering

\newcommand{\drawarm}[4]{
\node (#1m) at (#2, #3) [inner sep=2pt,circle,draw,fill]{};
\node (#1L) [above= #4 of #1m] {}; \draw[thick,-|] (#1m) -- (#1L);
\node (#1U) [below= #4 of #1m] {}; \draw[thick,-|] (#1m) -- (#1U);
}

\begin{tikzpicture}

\begin{scope}[red]
\drawarm{a11}{0}{0.95+.4}{1cm}
\end{scope}
\drawarm{a12}{1}{2.1+.4}{2cm}
\drawarm{a13}{2}{1.5+.4}{1.5cm}

\draw[dashed] (3,0) -- (3,7);

\drawarm{a21}{4}{3.5}{1cm}
\drawarm{a22}{5}{4.5}{2cm}
\begin{scope}[green]
\drawarm{a23}{6}{2.9}{.5cm}
\end{scope}

\draw[dashed] (7,0) -- (7,7);

\drawarm{a31}{8}{2.1+.1}{1cm}
\begin{scope}[red]
\drawarm{a32}{9}{1.5+.1}{.5cm}
\end{scope}
\drawarm{a33}{10}{3.05+.1}{2cm}

\draw[help lines] (-1, 2.3) -- (11, 2.3);
\draw[help lines] (-1, 2.5) -- (11, 2.5);
\node at (11.2,2.4) [anchor=center,font=\scriptsize] {$\}\epsilon$};

\end{tikzpicture}
\caption{Stopping rule \eqref{eq:stopping}. The algorithm stops because the lower bound of the green arm beats up to slack $\epsilon$ the upper bound for at least one arm (marked red) in each other action. In this case action $\ihat = 2$ is recommended.}\label{fig:stoppingrule}
\end{figure}

Concerning the stopping rule, which depends on the parameter $\epsilon \geq 0$ ($\epsilon$ can be set to zero if $\mu_{1,1}>\mu_{2,1}$), it is defined as the first moment when, according to the confidence intervals, some action $\ihat$ is probably approximately better than all other actions' best responses:
\begin{equation}\label{eq:stopping}
	\tau = \inf \left\{ t \in 2\N : \min_i \left[ \max_{i' \neq i} \min_{1 \leq j' \leq K_{i'}} \UCB_{i',j'}(t) - \min_{1 \leq j \leq K_i} \LCB_{i,j}(t)\right] < \epsilon\right\}\;.
\end{equation}
Then arm $\ihat = \ihat(\tau)$,  the empirical maximin action at that time, is recommended to player $A$. The stopping rule is illustrated in Figure~\ref{fig:stoppingrule}. With the notation of the sampling rule, this amounts to stopping when
$\LCB_{H_t}(t) > \UCB_{S_t}(t) - \epsilon$.

\subsection{Analysis of the Algorithm}

We analyze the algorithm under the assumptions $\mu_{1,1}<\mu_{2,1}$ and $\epsilon=0$. We consider the Hoeffding-type confidence bounds
\begin{equation}\LCB_{P}(t) = \hat{\mu}_P(t) - \sqrt{\frac{\beta(t,\delta)}{2N_P(t)}} \ \ \ \text{and} \ \ \ \UCB_P(t) = \hat{\mu}_P(t) + \sqrt{\frac{\beta(t,\delta)}{2N_P(t)}}\;,\label{bounds:Subgaussian}\end{equation}
where $\beta(t,\delta)$ is some exploration rate. A choice of $\beta(t,\delta)$ that ensures the $\delta$-PAC property~\eqref{eq:PAC} is given below. In order to highlight the dependency of the stopping rule on the risk level $\delta$, we denote it by $\tau_\delta$. 
\begin{theorem}\label{thm:subSC} Let
	\[H^*(\bm\mu) = \sum_{(i,j)\in \cP}\frac{1}{\max \left[\left(\mu_{i,1} - \frac{\mu_{1,1} + \mu_{2,2}}{2}\right)^2, (\mu_{i,j} - \mu_{i,1})^2\right]}.\]
	On the event 
	\[\cE = \bigcap_{P \in \cP} \bigcap_{t \in 2\N} \Big\{\mu_P \in [\LCB_P(t),\UCB_P(t)]\Big\}\;,\]
	the M-LUCB strategy returns the maximin action and uses a total number of samples upper-bounded by 
 \[T(\bm \mu,\delta) = \inf \big\{ t \in \N : 4 H^*(\bm\mu) \beta(t,\delta) < t \big\}.\]
\end{theorem}
According to~Theorem~\ref{thm:subSC}, the exploration rate should be large enough to control $\bP_{\bm\mu}(\cE)$, and as small as possible so as to minimize $T(\bm\mu,\delta)$.
The self-normalized deviation bound of~\cite{KLUCBJournal} gives a first solution (Corollary~\ref{cor:KLUCB}), whereas Lemma 7 of~\cite{JMLR15} yields Corollary~\ref{cor:NewSelf}. In both cases, explicit bounds on $T(\bm \mu,\delta)$ are obtained using the technical Lemma~\ref{lem:technical} stated in Appendix \ref{sec:proofsMLUCB}.

\begin{corollary}\label{cor:KLUCB} Let $\alpha>0$ and $C=C_\alpha$ be such that
	\[e \overline{K}\sum_{t=1}^\infty\frac{(\log t)(\log(C t^{1+\alpha}))}{t^{1+\alpha}} \leq C\;, \]
	and $\delta$ such that $4(1+\alpha) (C/\delta)^{1/(1+\alpha)} > 4.85$. 
	With probability larger than $1-\delta$, the M-LUCB strategy using the exploration rate 
	\begin{equation}\beta(t,\delta) = \log \left(\frac{Ct^{1+\alpha}}{\delta}\right)\;,\label{ExploBasic}\end{equation}
	 returns the maximin action within a number of steps upper-bounded as  
	\[\tau_\delta \leq 4 H^*(\bm\mu) \left[\log\left(\frac{1}{\delta}\right) +\log(C(4(1+\alpha)H^*(\bm\mu))^{1+\alpha}) + 2(1+\alpha)\log\log \left(\frac{4(1+\alpha)H^*(\bm\mu)C^{\frac{1}{1+\alpha}}}{\delta^{\frac{1}{1+\alpha}}}\right)\right]\] 
\end{corollary}

\begin{corollary}\label{cor:NewSelf} For $b,c$ such that $c>2$ and $b>c/2$, let the exploration rate be
	\[\beta(t,\delta) = \log \frac{1}{\delta} + b \log\log \frac{1}{\delta} + c \log\log(et)\]
	and 
	\[f_{b,c}(\delta) = \overline{K}\sqrt{e}\frac{\pi^2}{3}\frac{1}{8^{c/2}} \frac{(\sqrt{\log(1/\delta)+b\log\log(1/\delta)} + 2\sqrt{2})^c}{(\log(1/\delta))^b}\;,\]
	then with probability larger than $1- f_{b,c}(\delta)\delta$, M-LUCB returns the maximin action and, for some positive constant $C_c$ and for $\delta$ small enough, 
		\[\tau_\delta \leq 4 H^*(\bm\mu) \left[\log\left(\frac{1}{\delta}\right) +\log(8C_c H^*(\bm\mu)) + 2\log\log \left(\frac{8C_cH^*(\bm\mu)}{\delta}\right)\right]\] 
\end{corollary}
Elaborating on the same ideas, it is possible to obtain results in expectation, at the price of a less explicit bound, that holds for a slightly larger exploration rate. 

\begin{theorem}\label{thm:SCExpectation} The M-LUCB algorithm using $\beta(t,\delta)$ defined by \eqref{ExploBasic}, with $\alpha > 1$,  is $\delta$-PAC and satisfies  
\[\limsup_{\delta \rightarrow 0}\frac{\bE_{\bm \mu} [\tau_\delta]}{\log (1/\delta)} \leq 4 H^*(\bm\mu).\] 
\end{theorem}
The complexity term $H^*(\bm\mu)$ is easy to interpret: the number of draws of an arm $(i,j)$ is upper bounded by the typical number of samples needed to either discriminate $\mu_{i,j}$ from the smallest arm associated to the same action, $\mu_{i,1}$, or to discriminate $\mu_{i,1}$ from a `virtual arm' with mean $(\mu_{1,1} + \mu_{2,1})/2$. We view this virtual arm (that corresponds to the choice of a parameter $c$ in Appendix~\ref{sec:proofsMLUCB}) as an artifact of our proof, and we conjecture that it could be replaced by $\mu_{2,1}$ for arms in $\Arms_1$ and by $\mu_{1,1}$ for other arms. In the particular case of two actions by players, we propose the following finer result, that holds for the variant of M-LUCB that samples the least drawn arm among $H_t$ and $S_t$ at round $t+1$.

\begin{theorem}\label{thm:ParticularCase2} Assume  $K=K_1=K_2=2$. The M-LUCB algorithm using $\beta(t,\delta)$ defined by \eqref{ExploBasic} with $\alpha > 1$ is $\delta$-PAC and satisfies
\[\limsup_{\delta \rightarrow 0}\frac{\bE_{\bm\mu}[\tau_\delta]}{\log (1/\delta)} \leq 8\left[\frac{2}{(\mu_{1,1} - \mu_{2,1})^2} + \frac{1}{(\mu_{1,2} - \mu_{2,1})^2} + \frac{1}{\max \left[(\mu_{1,1} - \mu_{2,1})^2,(\mu_{2,2}-\mu_{2,1})^2\right]}\right].\] 
\end{theorem}

\subsection{Improved Intervals and Stopping Rule}\label{sec:NewTools}

The symmetry and the simple form of the sub-gaussian confidence intervals~\eqref{bounds:Subgaussian} are convenient for the analysis, but they can be greatly improved thanks to better deviation bounds for Bernoulli distributions. A simple improvement (see~\cite{COLT13}) is to use Chernoff confidence intervals, based on the binary relative entropy function $d(x,y)=x\log(x/y) + (1-x)\log((1-x)/(1-y))$. 
Moreover,  the use of a better stopping rule based on generalized likelihood ratio tests (GLRT) has been proposed recently for best-arm identification, leading to significant improvements. We propose here an adaptation of the Chernoff stopping rule of~\cite{GK16}, valid for the case $\epsilon=0$.  

This stopping rule based on the statistic: 
\[Z_{P,Q}(t) := \log\frac{\max_{\mu_P'\geq \mu_Q'} p_{\mu_P'}\left(\underline{X}^P_{N_P(t)}\right)p_{\mu_Q'}\left(\underline{X}^Q_{N_Q(t)}\right)}{\max_{\mu_P'\leq \mu_Q'} p_{\mu_P'}\left(\underline{X}^P_{N_P(t)}\right)p_{\mu_Q'}\left(\underline{X}^Q_{N_Q(t)}\right)}\;,\]
where $\underline{X}^P_s$ is a vector that contains the first $s$ observations of arm $P$ and $p_{\mu}(Z_1,\dots,Z_s)$ is the likelihood of $s$ i.i.d. observations from a Bernoulli distribution with mean $\mu$. Introducing the weighted sum of empirical means of two arms,  
\begin{eqnarray*}
	\hat{\mu}_{P,Q}(t) &:=& \frac{N_P(t)}{N_P(t) + N_Q(t)}\hat{\mu}_P(t) + \frac{N_Q(t)}{N_P(t) + N_Q(t)}{\hat{\mu}_Q(t)},
\end{eqnarray*}
it appears that for $\hat{\mu}_P(t) \geq \hat{\mu}_Q(t)$,  
\[Z_{P,Q}(t) = N_P(t) d\left(\hat{\mu}_P(t), \hat{\mu}_{P,Q}(t)\right) + N_Q(t) d\left(\hat{\mu}_Q(t), \hat{\mu}_{P,Q}(t) \right)\;,\]
and $Z_{P,Q}(t) = - Z_{Q,P}(t)$.
The stopping rule is defined as 
\begin{align}
	\tau & = \inf \Big\{t \in \N : \exists i \in \{1,\dots,K\} : \forall i' \neq i, \exists j'\in \{1,\dots,K_{i'}\} : \forall j \in \{1,\dots,K_i\}, Z_{(i,j), (i',j')}(t) > \beta(t,\delta)\Big\} \nonumber\\
	&= \inf \Big\{t \in \N : \max_{i \in \{1,\dots,K\}} \min_{i' \neq i} \max_{j'\in\{1,\dots,K_{i'}\}}\min_{j\in\{1,\dots,K_i\}}  Z_{(i,j), (i',j')}(t) > \beta(t,\delta)\Big\}\;.\label{StoppingChernoff}
\end{align}

\begin{proposition}\label{prop:PACbetter} Using the stopping rule~\eqref{StoppingChernoff} with the exploration rate $\beta(t,\delta) = \log \left(\frac{2 K_1 (K-1) t}{\delta}\right)$, whatever the sampling rule, if $\tau$ is a.s. finite, the recommendation is correct with probability $\bP_{\bm\mu}\left(\ihat = i^*\right) \geq 1- \delta.$
\end{proposition}

\paragraph{Sketch of Proof.}
Recall that in our notation the optimal action is $i^*=1$.
\begin{eqnarray*}
	\bP_{\bm\mu}\left(\ihat  \neq 1\right) & \leq & \bP_{\bm\mu}\left(\exists t\in \N, \exists i \in \{1,\dots,K\}\setminus \{1\}, \exists j \in \{1,\dots,K_1\},  Z_{(i,1),(1,j)}(t) > \beta(t,\delta) \right) \\
	& \leq & \sum_{i=2}^K \sum_{j=1}^{K_1} \bP_{\bm\mu}\left(\exists t\in \N, Z_{(i,1),(1,j)}(t) > \beta(t,\delta) \right)\;.
\end{eqnarray*}
Note that for $i\neq 1$, $\mu_{(i,1)}< \mu_{(1,j)}$ for all $j \in \{1,\dots,K_1\}$. The result follows from the following bound proved in \cite{GK16}: whenever $\mu_P< \mu_Q$, for any sampling strategy,
\begin{equation}\bP_{\bm\mu}\left(\exists t \in \N : Z_{P,Q}(t) > \log\left(\frac{2t}{\delta}\right)\right) \leq \delta\;.\label{ineq:ChernoffRule}\end{equation}

\section{A Racing algorithm}\label{sec:MRacing}
We now propose a Racing-type algorithm for the maximin action identification problem, inspired by another line of algorithms for best arm identification \citep{MaronMoore:97,EvenDaral06,COLT13}. Racing algorithms are simple and powerful methods that progressively concentrate on the best actions. We give in this section an analysis of a Maximin-Racing algorithm that relies on the refined information-theoretic tools introduced in the previous section.    


\subsection{A generic Maximin-Racing Algorithm}
The Maximin Racing algorithm maintains a set of active arms $\cR$ and proceeds in rounds, in which all the active arms are sampled. 
At the end of round $r$, all active arms have been sampled $r$ times and some arms may be eliminated according to some \emph{elimination rule}.
We denote by $\hat{\mu}_{P}(r)$ the average of the $r$ observations on arm $P$.
The elimination rule relies on an \emph{elimination function} $f(x,y)$ ($f(x,y)$ is large if $x$ is significantly larger than $y$), and on a \emph{threshold function} $\beta(r,\delta)$. 

The Maximin-Racing algorithm presented below performs two kinds of eliminations: the largest arm in each set $\cR_i$ may be eliminated if it appears to be significantly larger than the smallest arm in $\cR_i$ (\emph{high arm elimination}), and the group of arms $\cR_i$ containing the smallest arm may be eliminated (all the arms in $\cR_i$ are removed from the active set) if it contains one arm that appears significantly smaller than all the arms of another group $\cR_j$ (\emph{action elimination}).

\paragraph{Maximin Racing algorithm} 
~\\\emph{Parameters.} Elimination function $f$, threshold  function $\beta$
\\\emph{Initialization.} For each $i\in\{1,\dots,K\},\cR_i = \Arms_i$, and  $\cR := \cR_1 \cup \dots \cup \cR_K$.
\\\emph{Main Loop.} At round $r$:
\begin{itemize}
 \item all arms in $\cR$ are drawn, empirical means $\hat{\mu}_{P}(r)$, $P\in \cR$ are updated
 \item \emph{High arms elimination step}: 
 for each action $i=1\dots K$, if $|\cR_i| \geq 2$ and  
 \begin{equation}r f \left(\max_{P \in \cR_i} \hat{\mu}_{P}(r),\min_{P \in \cR_i} \hat{\mu}_{P}(r)\right) \geq \beta(r,\delta)\;,\label{LargeArmsElimination}\end{equation}
 then remove $P_m = \argmax{j \in \cR_i} \ \hat{\mu}_{P}(r)$ from the active set : $\cR_i = \cR_i \backslash \{P_m\}$, $\cR = \cR \backslash \{P_m\}$.
 
  \item \emph{Action elimination step}: 
  if $(\itilde,\jtilde) = \argmin{P \in \cR} \ \hat{\mu}_{P}(r)$ and if
  \[r f \left(\max_{i\neq \itilde} \min_{P \in \cR_i} \hat{\mu}_{P}(r), \hat{\mu}_{(\itilde,\jtilde)}(r)\right) \geq \beta(r,\delta)\;,\]
 then remove $\itilde$ from the possible maximin actions: $\cR = \cR \backslash \cR_{\itilde}$ and $\cR_{\itilde} = \emptyset$.
\end{itemize}
The algorithm stops when all but one of the $\cR_i$ are empty, and outputs the index of the remaining set as the maximin action. If the stopping condition is not met for 
\[r = r_0: = \frac{2}{\epsilon^2}\log \left(\frac{4\overline{K}}{\delta}\right),\]
then the algorithm stops and returns one of the empirical maximin actions.

\subsection{Tuning the Elimination and Threshold Functions} 

In the best-arm identification literature, several elimination functions have been studied. The first idea, presented in the Successive Elimination algorithm of \cite{EvenDaral06}, is to use the simple difference $f(x,y)=(x-y)^2\ind_{(x \geq y)}$; in order to take into account possible differences in the deviations of the arms, the KL-Racing algorithm of \cite{COLT13} uses an elimination function equivalent to
$f(x,y)=d_*(x,y)\ind_{(x \geq y)}$, where $d_*(x,y)$ is defined as the common value of $d(x,z)$ and $d(y,z)$ for the unique $z$ satisfying $d(x,z)=d(y,z)$. In this paper, we use the divergence function
\begin{equation}f(x,y) = I(x,y) : = \left[d\left(x,\frac{x+y}{2}\right) + d\left(y,\frac{x+y}{2}\right)\right]\ind_{(x \geq y)}\label{def:EliminationRule}\end{equation}
inspired by the deviation bounds of Section~\ref{sec:NewTools}. In particular, using again Inequality~\eqref{ineq:ChernoffRule} for the uniform sampling rule yields, whenever $\mu_P < \mu_Q$, 
\begin{equation}\bP_{\bm\mu}\left(\exists r \in \N :  r I( \hat{\mu}_{P}(r), \hat{\mu}_{Q}(r)) \geq \log \frac{2r}{\delta}\right) \leq \delta.\label{inequ:BBlock}\end{equation}
Using this bound, Proposition~\ref{lem:PACRacing} (proved in Appendix~\ref{proof:PACRacing}) proposes a choice of the threshold function for which the Maximin-Racing algorithm is $\delta$-PAC.   
\begin{proposition}\label{lem:PACRacing} With the elimination function $I(x,y)$ of Equation~\eqref{def:EliminationRule} and with the threshold function $\beta(t,\delta) = \log \left(4{C_K} t /\delta\right)$, the Maximin-Racing algorithm satisfies 
\[\bP_{\bm\mu} \left(\mu_{1,1} - \mu_{\ihat,1} \leq \epsilon \right) \geq 1- \delta,\]
with $C_K \leq (\overline{K})^2$. If $\mu_{1,1}>\mu_{1,2}$ and if $\forall i, \mu_{i,1} < \mu_{i,2}$, then $C_K = K \times \max_i K_i$. 
\end{proposition}

\subsection{Sample Complexity Analysis}

We propose here an asymptotic analysis of the number of draws of each arm $(i,j)$ under the Maximin-Racing algorithm, denoted by $\tau_{\delta}(i,j)$. These bounds are expressed with the deviation function $I$, and hold for $\epsilon>0$. For $\epsilon=0$, one can provide similar bounds under the additional assumption that all arms are pairwise distinct. 

\begin{theorem}\label{thm:PACRacing} Assume $\mu_{1,1} > \mu_{2,1}$. For every $\epsilon>0$, and for  $\beta(t,\delta)$ chosen as in Proposition~\ref{lem:PACRacing}, the Maximin-Racing algorithm satisfies 
\[\limsup_{\delta \rightarrow 0} \frac{\bE_{\bm\mu}[\tau_{\delta}(1,1)]}{\log(1/\delta)} \leq \frac{1}{\max\big(\epsilon^2/2, I(\mu_{2,1},\mu_{1,1})\big)}\]
and, for any $(i,j) \neq (1,1)$, 
\[\limsup_{\delta \rightarrow 0} \frac{\bE_{\bm\mu}[\tau_{\delta}(i,j)]}{\log(1/\delta)} \leq \frac{1}{\max\big(\epsilon^2/2, I(\mu_{i,1},\mu_{1,1}), I(\mu_{i,j},\mu_{i,1})\big)}.\]
\end{theorem}
It follows from Pinsker's inequality that $I(x,y) > (x-y)^2$, and hence Theorem~\ref{thm:PACRacing} implies in particular that for the M-Racing algorithm (for a sufficiently small $\epsilon$)
\begin{eqnarray*}\limsup_{\delta \rightarrow 0} \frac{\bE_{\bm\mu}[\tau_\delta]}{\log(1/\delta)} &\leq& \frac{1}{(\mu_{1,1} - \mu_{2,1})^2} + \sum_{j=2}^{K_1} \frac{1}{(\mu_{1,j} - \mu_{1,1})^2} +
\sum_{i=2}^{K}\sum_{j=1}^{K_i} \frac{1}{(\mu_{1,1} - \mu_{i,1})^2 \vee (\mu_{i,j} - \mu_{i,1})^2}.
\end{eqnarray*}
The complexity term on the right-hand side is reminiscent of the quantity $H^*(\bm\mu)$ introduced in Theorem~\ref{thm:subSC}. The terms corresponding to arm in $\cP \setminus \Arms_1$ are comparable to the corresponding terms in $H^*(\bm\mu)$ (they are actually strictly smaller since no `virtual arm' $(\mu_{1,1} + \mu_{2,1})/2$ have been introduced in the analysis of M-Racing). However, the terms corresponding to the arms $(1,j),j\geq 2$ are strictly larger than the corresponding terms in $H^*(\bm\mu)$.  But this is mitigated by the fact that there is no multiplicative constant in front of the complexity term. Besides, as Theorem~\ref{thm:PACRacing} involves the deviation function $I(x,y) = d(x,(x+y)/2) + d(y,(x+y)/2$ and not a subgaussian approximation, they can indeed be significantly better.

\section{Numerical Experiments and Discussion}\label{sec:Experiments}
In the previous sections, we have proposed two different algorithms for the maximin action identification problem. The analysis that we have given does not clearly advocate the superiority of one or the other.
The goal of this section is to propose a brief numerical comparison in different settings, and to compare with other possible strategies.

We will notably study empirically two interesting variants of M-LUCB. The first improvement that we propose is the  M-KL-LUCB strategy, based on  KL-based  confidence bounds (\cite{COLT13}). The second variant, M-Chernoff, additionally improves the stopping rule as presented in Section~\ref{sec:NewTools}. Whereas Proposition~\ref{prop:PACbetter} justifies the use of the exploration rate $\beta(t,\delta)=\log(4\overline{K}^2 t / \delta)$, which is over-conservative in practice, we use  $\beta(t,\delta)=\log((\log(t)+1)/\delta)$ in all our experiments, as suggested by Corollary~\ref{cor:NewSelf} (this appears to be already quite a conservative choice in practice). In the experiments, we set $\delta=0.1$, $\epsilon=0$.

To simplify the discussion and the comparison, we first focus on the particular case in which there are two actions for each player.  
As an element of comparison, one can observe that finding $i^*$ is at most as hard as finding the worst arm (or the three best) among the four arms $(\mu_{i,j})_{1\leq i,j \leq 2}$. Thus, one could use standard best-arm identification strategies like the (original) LUCB algorithm. For the latter, the complexity is of order
\[\frac{2}{(\mu_{1,1} - \mu_{2,1})^2} + \frac{1}{(\mu_{1,2} - \mu_{2,1})^2} + \frac{1}{(\mu_{2,2} - \mu_{2,1})^2}\;,\]
which is much worse than the complexity term obtained for M-LUCB in Theorem~\ref{thm:ParticularCase2} when $\mu_{2,2}$ and $\mu_{2,1}$ are close to one another. This is because a best arm identification algorithm does not only find the maximin action, but additionally figures out which of the arms in the other action is worst. Our algorithm does not need to discriminate between $\mu_{2,1}$ and $\mu_{2,2}$, it only tries to assess that one of these two arms is smaller than $\mu_{1,1}$. However, for specific instances in which the gap between $\mu_{2,2}$ and $\mu_{2,1}$ is very large, the difference vanishes.
This is illustrated in the numerical experiments of Table~\ref{fig:2arms}, which involve the following three sets of parameters (the entry $(i,j)$ in each matrix is the mean $\mu_{i,j}$):
\[
\bm\mu_1 =   \left[\begin{array}{cc}
  0.4  & 0.5 \\      
  0.3 & 0.35 \\      
   \end{array} \right] \ \ \ 
\bm\mu_2 =   \left[
\begin{array}{cc}
  0.4  & 0.5 \\       0.3 & 0.45 \\      
\end{array}\right] \ \ \ 
      \bm\mu_3 =   \left[
\begin{array}{cc}
  0.4  & 0.5 \\       
  0.3 & 0.6 \\      
\end{array} \right]
\]

\begin{table}[h]
\begin{tabular}{|c||c|c|c|c||c|c|c|c||c|c|c|c|}
 \hline
 & $\tau_{1,1}$ & $\tau_{1,2}$ & $\tau_{2,1}$ & $\tau_{2,2}$ & $\tau_{1,1}$ & $\tau_{1,2}$ & $\tau_{2,1}$ & $\tau_{2,2}$& $\tau_{1,1}$ & $\tau_{1,2}$ & $\tau_{2,1}$ & $\tau_{2,2}$\\ 
 \hline
 M-LUCB & 1762 & 198 & 1761& 462 & 1761 & 197 & 1760 & 110 & 1755 & 197 & 1755 & 36 \\
 \hline 
 M-KL-LUCB & 762 & 92 & 733 & 237 & 743 & 92 & 743 & 54 & 735 & 93 & 740 & 16 \\
 \hline
 M-Chernoff & 315 & 59 & 291 & 136 & 325 & 61 & 327 & 41 & 321 & 61 & 326 & 13 \\
 \hline 
 M-Racing & 324 & 152 & 301 & 298 & 329 & 161 & 318 & 137 & 322 & 159 & 323 & 35 \\
 \hline 
 KL-LUCB & 351 & 64 & 3074 & 2768 & 627 & 83 & 841 & 187 & 684 & 88 & 774 & 32 \\
 \hline
\end{tabular}
\caption{\label{fig:2arms} Number of draws of the different arms under the models parameterized by $\bm \mu_1, \bm \mu_2, \bm \mu_3$ (from left to right), averaged over $N=10000$ repetitions}
\end{table}
We also perform experiments in a model with 3x3-actions with parameters: 
\[
\bm\mu =   \left[\begin{array}{ccc}
  0.45  & 0.5 & 0.55 \\      
  0.35 & 0.4 & 0.6 \\
  0.3 & 0.47 & 0.52 \\
   \end{array} \right] 
\]
Figure~\ref{fig:3Actions} shows that the  best three algorithms in the previous experiments behave as expected: the number of draws of the arms are ordered exactly as suggested by the bounds given in the analysis. 
\begin{figure}[h]
\[
 \tau_{\text{M-KLLUCB}} =\left[\begin{array}{ccc}
798 & 212 &  92 \\
752 & 248 &  22 \\
210 & 44 &  21 \\
\end{array} \right] \ \ \ \
 \tau_{\text{M-Ch.}} =   \left[\begin{array}{ccc}
367 & 131 & 67 \\
333  & 156 & 18 \\
129 &  31 & 17 \\ 
\end{array} \right] \ \ \ 
\tau_{\text{M-Racing}} =   \left[\begin{array}{ccc}
 472 &  291 & 173 \\
 337 & 337 & 42 \\
 161 & 185 & 71 \\
 \end{array}\right] 
 \]  
 \caption{\label{fig:3Actions} Number of draws of each arm under the bandit model $\bm\mu$, averaged of $N=10000$ repetitions}
\end{figure}
These experiments tend to show that, in practice, the best two algorithms are M-Racing and M-Chernoff, with a slight advantage for the latter. However, we did not provide theoretical sample complexity bounds for M-Chernoff, and it is to be noted that the use of Hoeffding bounds in the M-LUCB algorithm (that has been analyzed) is a cause of sub-optimality. Among the algorithms for which we provide theoretical sample complexity guarantees, the M-Racing algorithm appears to perform best. 

\section{Perspectives} \label{sec:Perspectives} 
To finish, let us sketch the (still speculative) perspective of an important improvement.
For simplicity, we focus on the case where each player chooses among only two possible actions, and we change our notation, using:  $\mu_1:=\mu_{1,1},\mu_2:=\mu_{1,2},\mu_3:=\mu_{2,1},\mu_4:=\mu_{2,2}$. 
As we will see below, the optimal strategy is going to depend a lot on the position of $\mu_4$ relatively to $\mu_1$ and $\mu_2$. 
Given $w = (w_1,\dots,w_4) \in \Sigma_K=\{w\in\R_+^4: w_1+\dots+w_4=1\}$, we define for  $a,b,c$ in $\{1,\dots,4\}$:
\[\mu_{a,b}(w) = \frac{w_a \mu_a + w_b \mu_b}{w_a + w_b} \ \ \text{and} \ \ \ 
\mu_{a,b,c}(w) = \frac{w_a \mu_a + w_b \mu_b + w_c \mu_c}{w_a + w_b + w_c}\;.
\]

Using a similar argument than the one of \cite{GK16} in the context of best-arm identification, one can prove the following (non explicit) lower bound on the sample complexity.

\begin{theorem} Any $\delta$-PAC algorithm satisfies 
\[\bE_{\bm\mu}[\tau_\delta] \geq T^*(\bm \mu)\, d(\delta,1-\delta),\]
where 
\begin{eqnarray}
 T^*(\bm \mu)^{-1} &: =& \sup_{w \in \Sigma_K} \inf_{\bm \mu' : \mu_1'\wedge \mu_2' < \mu_3' \wedge \mu_4'} \left(\sum_{a=1}^K w_a\, d(\mu_a,\mu_a')\right) \nonumber\\
&= &\sup_{w \in \Sigma_K} \min [F_1(\bm \mu,w), F_2(\bm \mu,w)],\label{GeneOptimPb}
\end{eqnarray}
where 
\[F_a(\bm \mu, w) = \left\{
\begin{array}{ll}
w_a\,d\big(\mu_a,\mu_{a,3}(w)\big) + w_3\,d\big(\mu_3,\mu_{a,3}(w)\big) & \text{if} \ \mu_4 \geq \mu_{a3}(w)\;, \\
w_a\,d\big(\mu_a,\mu_{a,3,4}(w)\big) + w_3\,d\big(\mu_3,\mu_{a,3,4}(w)\big) + w_4\,d\big(\mu_4,\mu_{a,3,4}(w)\big) &\text{otherwise}.
\end{array}
\right.\] 
\end{theorem}

\paragraph{A particular case.} When $\mu_4 > \mu_2$, for any $w\in\Sigma_K$ it holds that $\mu_4 \geq \mu_{1,3}(w)$ and $\mu_4 \geq \mu_{2,3}(w)$. Hence the complexity term can be rewritten to
 \[T^*(\bm \mu)^{-1} = \sup_{w \in \Sigma_K} \min_{a=1,2} w_a\,d\big(\mu_a,\mu_{a,3}(w)\big) + w_3\,d\big(\mu_3,\mu_{a,3}(w)\big)\;.
 \]
In that case it is possible to show that the following quantity, 
\[w^*(\bm \mu) = \argmax{w \in \Sigma_K} \ \min_{a=1,2}  w_a\,d\big(\mu_a,\mu_{a,3}(w)\,\big) + w_3\,d\big(\mu_3,\mu_{a,3}(w)\big)\]
is unique and to give a more explicit expression. This quantity is to be interpreted as the vector of proportions of draws of the arms by a strategy matching the lower bound.  In this particular case, one finds $w_4^*(\bm \mu) = 0$, showing that an optimal strategy could draw arm 4 only an asymptotically vanishing proportion of times as $\delta$ and $\epsilon$ go to~$0$.

\paragraph{Towards an Asymptotically Optimal Algorithm.} 

Assume that the solution of the general optimization  problem \eqref{GeneOptimPb} is well-behaved (unicity of the solution, continuity in the parameters,...) and that we can find an efficient algorithm to compute 
\[w^*(\bm \mu) = \argmax{w \in \Sigma_K} \min [F_1(\bm\mu,w), F_2(\bm\mu,w)]\]
for any given $\bm \mu$. In particular, for a fixed $w$ and $\bm \mu$, we need to be able to compute 
\[F(w,\bm \mu) = \inf_{\mu' \in \mathrm{Alt}(\bm \mu)} \sum_{a=1}^4 w_a d(\mu_a,\mu_a'),\]
where $\mathrm{Alt}(\bm \mu) = \{ \bm\mu' : i^*(\bm\mu) \neq i^*(\bm\mu')\}$. 
Then, if we can design a sampling rule ensuring that for all $a$, $N_a(t)/t$ tends to $w^*_a(\bm \mu)$, and if we combine it with the stopping rule 
\[\tau_\delta = \inf \bigg\{t \in \N : F\Big(\big(N_a(t)\big)_{a=1\dots 4} ,\, \hat{\bm\mu}(t)\Big) > \log(Ct / \delta)\bigg\}\]
for some positive constant $C$, then one could expect the following asymptotic optimality property: 
\[\limsup_{\delta\to 0}\frac{\bE_{\bm \mu}[\tau_\delta]}{\log(1/\delta)} \leq T^*(\bm \mu).\]
But proving that this stopping rule does ensures a $\delta$-PAC algorithm is not straightforward, and the analysis remains to be done.

\acks{
	This work was partially supported by the CIMI (Centre International de Math\'ematiques et d'{In\-for\-ma\-tique}) Excellence program while Emilie Kaufmann visited Toulouse in November 2015.
	The authors acknowledge the support of the French Agence Nationale de la Recherche (ANR), under grants ANR-13-BS01-0005 (project SPADRO) and ANR-13-CORD-0020 (project ALICIA). 
}

\bibliography{biblioBandits}

\appendix 

\section{Analysis of the Maximin-LUCB algorithm} \label{sec:proofsMLUCB}

We define the event 
\[\cE_t =\bigcap_{P\in \cP}(\mu_P \in [\LCB_P(t),\UCB_P(t)]),\]
so that the event $\cE$ defined in Theorem~\ref{thm:subSC} rewrites $\cE = \bigcap_{t\in 2\N} \cE_t$. 

Assume that the event $\cE$ holds. The arm $\ihat$ recommended satisfies, by definition of the algorithm, for all $i \neq \ihat$
\[\min_{j \in K_{\ihat}} \LCB_{(\ihat,j)}(\tau_\delta) > \min_{j \in K_i} \UCB_{(i,j)}(\tau_\delta) - \epsilon.\]
Using that $ \LCB_{P}(\tau_\delta) \leq \mu_P \leq \UCB_{P}(\tau_\delta)$ for all $P\in \cP$ (by definition of $\cE$) yields for all $i$
\[\mu_{\ihat,1}=\min_{j \in K_{\ihat}} \mu_{\ihat,j} > \min_{j \in K_i} \mu_{i,j} - \epsilon =  \mu_{i,1} - \epsilon,\]
hence $\max_{i \neq \ihat} \mu_{i,1} - \mu_{\ihat,1} < \epsilon$. Thus, either $\ihat = 1$ or $\ihat$ satisfies 
$ \mu_{1,1} - \mu_{\ihat,1} < \epsilon$. In both case, $\ihat$ is $\epsilon$-optimal, which proves that M-LUCB is correct on $\cE$. 

Now we analyze M-LUCB with $\epsilon=0$. Our analysis is based on the following two key lemmas, whose proof is given below. 

\begin{lemma}\label{lem:CoreLemma1} Let $c \in [\mu_{2,1},\mu_{1,1}]$ and $t\in 2\N$.  On $\cE_t$, if $(\tau_\delta > t)$, there exists $P \in \{ H_t,S_t\}$ such that 
	\[\left(c \in [\LCB_P(t), \UCB_P(t)]\right).\]
\end{lemma}

\begin{lemma}\label{lem:CoreLemma2} Let $c \in [\mu_{2,1},\mu_{1,1}]$ and $t\in 2\N$. On $\cE_t$, for every $(i,j) \in \{H_t,S_t\}$, 
	\[c \in [\LCB_{(i,j)}(t), \UCB_{(i,j)}(t)] \ \ \Rightarrow \ \ N_{(i,j)}(t) \leq \min\left(\frac{2}{(\mu_{i,1} - c)^2}, \frac{2}{(\mu_{i,j} - \mu_{i,1})^2}\right) \beta(t,\delta)\] 
\end{lemma}

Defining, for every arm $P\in\cP$ the constant 
\[c_P = \frac{1}{\max\left[\left(\mu_{i,1} - \frac{\mu_{1,1} + \mu_{2,1}}{2}\right)^2,(\mu_{i,j} - \mu_{i,1})^2\right]},\]
combining the two lemmas (for the particular choice $c=\frac{\mu_{1,1} + \mu_{2,1}}{2}$) yields the following key statement:
\begin{equation} \label{KeyStatement}
\cE_t \cap (\tau_\delta > t) \ \ \Rightarrow \ \exists P \in \{H_t,S_t\} : N_P(t) \leq 2 c_P \beta(t,\delta).\end{equation}
Note that $H^*(\bm\mu) = \sum_{P\in\cP} c_P$, from its definition in Theorem~\ref{thm:subSC}.

\subsection{Proof of Theorem~\ref{thm:subSC}}

Let $T$ be a deterministic time. On the event $\cE = \bigcap_{t\in 2\N}\cE_t$, using \eqref{KeyStatement} and the fact that for every even $t$, $(\tau_\delta > t)=(\tau_\delta > t+1)$ by definition of the algorithm, one has
\begin{eqnarray*} 
	\min(\tau_\delta,T)  &=& \sum_{t=1}^T \ind_{(\tau_\delta > t)} =2\sum_{\substack{t \in 2\N\\ t \leq T}} \ind_{(\tau_\delta > t)} =2 \sum_{\substack{t \in 2\N\\ t \leq T}} \ind_{\left(\exists P \in \{H_t,S_t\} : N_P(t) \leq 2c_p \beta(t,\delta)\right)} \\
	& \leq & 2 \sum_{\substack{t \in 2\N\\ t \leq T}} \sum_{P \in \cP} \ind_{(P_{t+1}=P)\cup(P_{t+2}=P)}\ind_{(N_P(t) \leq 2 c_P \beta(T,\delta))} \\ 
	& \leq & 4\sum_{P \in \cP} c_P \beta(T,\delta) = 4 H^*(\bm\mu) \beta(T,\delta). 
\end{eqnarray*}
For any $T$ such that 
$4H^*(\bm\mu) \beta(T,\delta) < T$, one has $\min(\tau_\delta, T)< T$, which implies $\tau_\delta <  T$. Therefore $\tau_\delta \leq T(\bm\mu,\delta)$ for $T(\bm\mu,\delta)$ defined in Theorem~\ref{thm:subSC}.

\subsection{Proof of Theorem~\ref{thm:SCExpectation}}

Let $\gamma > 0$. Let $T$ be a deterministic time. On the event $\cG_T = \bigcap_{\substack{t \in 2\N \\ \lfloor \gamma T \rfloor \leq t \leq T}} \cE_t$, one can write 

\begin{eqnarray*} 
	\min(\tau_\delta,T)  &=& 2\gamma T + 2\sum_{\substack{t \in 2\N \\ \lfloor \gamma T \rfloor\leq t \leq T}} \ind_{(\tau_\delta > t)} =2\gamma T + 2\sum_{\substack{t \in 2\N \\ \lfloor \gamma T \rfloor\leq t \leq T}} \ind_{\left(\exists P \in \{H_t,S_t\} : N_P(t) \leq 2c_p \beta(t,\delta)\right)} \\
	& \leq & 2\gamma T + 2\sum_{\substack{t \in 2\N \\ \lfloor \gamma T \rfloor\leq t \leq T}}\sum_{P \in\cP} \ind_{(P_{t+1}=P)\cup(P_{t+1}=P)}\ind_{(N_P(t) \leq 2 c_P \beta(T,\delta))} \\ 
	& \leq & 2\gamma T + 4 H^*(\bm\mu) \beta(T,\delta). 
\end{eqnarray*}
Introducing $T_\gamma(\bm\mu,\delta) := \inf\{ T \in \N : 4H^*(\bm\mu) \beta(T,\delta) < (1-2\gamma)T\}$, for all $T \geq T_\gamma(\bm\mu,\delta)$, $\cG_T \subseteq (\tau_\delta \leq T)$.  One can bound the expectation of $\tau_\delta$ in the following way (using notably the self-normalized deviation inequality of \cite{KLUCBJournal}): 

\begin{eqnarray*} \bE_{\bm\mu}[\tau_\delta] &= &\sum_{T=1}^\infty \bP_{\bm\mu}(\tau_\delta > T) \leq T_\gamma + \sum_{T=T_\gamma}^\infty \bP_{\bm\mu}\left(\tau_\delta > T\right) \leq   T_\gamma + \sum_{T=T_\gamma}^\infty \bP_{\bm\mu}\left(\cG_T^c\right)\\
	& \leq & T_\gamma + \sum_{T=1}^\infty \sum_{t=\gamma T}^T \sum_{P\in\cP} \left[\bP_{\bm\mu}\left(\mu_P > \hat{\mu}_P(t) + \sqrt{\frac{\beta(t,\delta)}{2N_P(t)}}\right)+\bP_{\bm\mu}\left(\mu_P < \hat{\mu}_P(t) -  \sqrt{\frac{\beta(t,\delta)}{2N_P(t)}}\right)\right] \\
	& \leq & T_\gamma + \sum_{T=1}^\infty \sum_{t=\gamma T}^T 2 \overline{K} \bP_{\bm\mu}\left(\mu_P > \hat{\mu}_P(t) + \sqrt{\frac{\beta(t,1)}{2N_P(t)}}\right) \\
	& \leq & T_\gamma + \sum_{T=1}^\infty \sum_{t=\gamma T}^T 2\overline{K} e\log(t) \beta(t,1)\exp(-\beta(t,1)) \\
	& \leq & T_\gamma + \sum_{T=1}^\infty 2\overline{K} e T \log(T) \beta(T,1)\exp(-\beta(\gamma T ,1)) \\
	& = &  T_\gamma + \sum_{T=1}^\infty \frac{2\overline{K} e T \log(T) \log(CT^{1+\alpha})}{C\gamma^{1+\alpha}T^{1+\alpha}}, 
\end{eqnarray*}
where the series is convergent for $\alpha > 1$. One has 
\[T_\gamma(\bm \mu,\delta) = \inf\left\{ T \in \N : \log\left(\frac{CT^{1+\alpha}}{\delta}\right) < \frac{(1-2\gamma)T}{4H^*(\bm\mu) }\right\}.\]
The technical Lemma~\ref{lem:technical} below permits to give an upper bound on $T_\gamma(\bm \mu,\delta)$ for small values of $\delta$, that implies in particular 
\[\limsup_{\delta \rightarrow 0} \frac{\bE_{\bm\mu}[\tau_\delta]}{\log(1/\delta)} \leq \frac{4 H^*(\bm\mu)}{1-2\gamma}.\]
Letting $\gamma$ go to zero yields the result.

\begin{lemma}\label{lem:technical} 
If  $\alpha, c_1,c_2>0$ are such that $a=(1+\alpha)c_2^{1/(1+\alpha)}/c_1>4.85$, then
\[x = \frac{1+\alpha}{c_1}\big(\log(a)+2\log(\log(a))\big)\]
is such that $c_1  x \geq \log(c_2 x^{1+\alpha})$. 
\end{lemma}
\paragraph{Proof.} One can check that if $a\geq 4.85$, then $\log^2(a)>\log(a)+2\log(\log(a))$. 
Thus, $y=\log(a)+2\log(\log(a))$ is such that $y\geq \log(ay)$.
Using $y=c_1x/(1+\alpha)$ and $a=(1+\alpha)c_2^{1/(1+\alpha)}/c_1$, one obtains the result.
\qed

\subsection{Proof of Lemma~\ref{lem:CoreLemma1}}

We show that on $\cE_t \cap (\tau_\delta > t)$, the following four statements cannot occur, which yields that the threshold $c$ is contained in one of the intervals $\cI_{H_t}(t)$ or $\cI_{S_t}(t)$: 
\begin{enumerate}
	\item $(\LCB_{H_t}(t) > c)\cap (\LCB_{S_t}(t) > c)$ 
	\item $(\UCB_{H_t}(t) < c)\cap (\UCB_{S_t}(t) < c)$ 
	\item $(\UCB_{H_t}(t) < c)\cap (\LCB_{S_t}(t) > c)$ 
	\item $(\LCB_{H_t}(t) > c)\cap (\UCB_{S_t}(t) < c)$ 
\end{enumerate}

1. implies that there exists two actions $i$ and $i'$ such that $\forall j \leq K_i, \LCB_{i,j}(t) \geq c $ and $\forall j' \leq K_{i'}, \LCB_{i',j'}(t) \geq c$. Because $\cE_t$ holds, one has in particular $\mu_{i,1}> c$ and $\mu_{j,1} > c$, which is excluded since $\mu_{1,1}$ is the only such arm that is larger than $c$. 

2. implies that for all $i \in \{1,K\}$, $\UCB_{(i,c_i(t))}(t) \leq c$. Thus, in particular $\UCB_{(1,c_1(t))} \leq c$ and, as $\cE_t$ holds, there exists $j \leq K_1$ such that $\mu_{1,j} < c$, which is excluded. 

3. implies that there exists $i \neq \ihat(t)$ such that $\min_{j} \hat{\mu}_{i,j}(t) > \hat{\mu}_{H_t}(t) \geq  \min_{j}\hat{\mu}_{(\ihat(t),j)}(t)$, which contradicts the definition of $\ihat(t)$. 

4. implies that $\UCB_{H_t}(t) > \LCB_{S_t}(t)$, thus the algorithm must have stopped before the  $t$-th round, which is excluded since $\tau_\delta > t$.

We proved that there exists $P \in \{H_t,S_t\}$ such that $c\in \cI_{P}(t)$. 


\subsection{Proof of Lemma~\ref{lem:CoreLemma2}} 

Assume that $\cE_t$ holds and that $c \in [\LCB_{(i,j)}(t), \UCB_{(i,j)}(t)]$. We first show that $(i,1)$ is also contained in  $[\LCB_{(i,j)}(t), \UCB_{(i,j)}(t)]$. First, by definition of the algorithm, if $(i,j)=H_t$ or $S_t$, one has $(i,j)=(i,c_i(t))$, hence 
\[\LCB_{(i,j)}(t) \leq \LCB_{(i,1)}(t) \leq \mu_{i,1},\]
using that $\cE_t$ holds. Now, if we assume that $\mu_{i,1} > \UCB_{(i,j)}(t)$, because $\cE_t$ holds, one has $\mu_{i,1}>\mu_{i,j}$, which is a contradiction. Thus, $\mu_{i,1} \leq \UCB_{(i,j)}(t)$. 

As $c$ and $\mu_{i,1}$ are both contained in $[\LCB_{(i,j)}(t), \UCB_{(i,j)}(t)],$ whose diameter is $2\sqrt{\beta(t,\delta)/(2N_{(i,j)}(t))}$, one has 
\[|c - \mu_{i,1}| < 2\sqrt{\frac{\beta(t,\delta)}{2N_{(i,j)}(t)}} \ \ \Leftrightarrow \ \ \ N_{(i,j)}(t) \leq \frac{2\beta(t,\delta)}{(\mu_{i,1}-c)^2}.\]
Moreover, one can use again that $\LCB_{(i,j)}(t) \leq \LCB_{(i,1)}(t)$ to write
\begin{eqnarray*}
	\UCB_{(i,j)}(t)  - 2\sqrt{\frac{\beta(t,\delta)}{2N_{(i,j)}(t)}}  & \leq & \LCB_{(i,1)}(t) \\
	\mu_{i,j} - 2\sqrt{\frac{\beta(t,\delta)}{2N_{(i,j)}(t)}} & \leq & \mu_{i,1},
\end{eqnarray*}
which yields $N_{(i,j)}(t) \leq \frac{2\beta(t,\delta)}{(\mu_{i,j}-\mu_{i,1})^2}$ and concludes the proof. 

\subsection{Proof of Theorem~\ref{thm:ParticularCase2}}

In the particular case of two actions by player, we analyze the version of LUCB that draws only one arm per round. More precisely, in this particular case, letting 
\[X_t = \argmin{j=1,2} \ \LCB_{(1,j)}(t) \ \ \ \text{and} \ \ \ Y_t = \argmin{j=1,2} \ \LCB_{(2,j)}(t),\]
one has $P_{t+1} = \argmax{P \in \{X_t,Y_t\}} \ N_P(t)$.

The analysis follows the same lines as that of Theorem~\ref{thm:SCExpectation}. First, we notice that the algorithm outputs the maximin action on the event $\cE = \cap_{t\in\N} \cE_t$, and thus the exploration rate defined in Corollary~\ref{cor:KLUCB} guarantees a $\delta$-PAC algorithm. Then, the sample complexity analysis relies on a specific characterization of the draw of each of the arms given in Lemma~\ref{lem:Events} below (which is a counterpart of Lemma~\ref{lem:CoreLemma2}). This result justifies the new complexity term that appears in Theorem~\ref{thm:ParticularCase2}. 

\begin{lemma}\label{lem:Events} On the event $\cE$, for all $P \in \cP$, one has 
\[(P_{t+1} = P)\cap (\tau_\delta > t) \subseteq \left(N_P(t) \leq 8 c_P \beta(t,\delta)\right),\]
with
\[c_{(1,1)} =  \frac{1}{(\mu_{1,1} - \mu_{2,1})^2}, \ \ \ c_{(1,2)} = \frac{1}{(\mu_{1,2} - \mu_{2,1})^2}, \ \ \ c_{(2,1)} =  \frac{1}{(\mu_{1,1} - \mu_{2,1})^2},\]
and
\[ c_{(2,2)} = \frac{1}{\min(4(\mu_{2,2} - \mu_{2,1})^2,(\mu_{1,1} - \mu_{2,1})^2)}.\]
\end{lemma}

\paragraph{Proof of Lemma~\ref{lem:Events}.} The proof of this result uses extensively the fact that the confidence intervals in \eqref{bounds:Subgaussian} are symmetric:
\[\UCB_P(t) = \LCB_P(t) + 2\sqrt{\frac{\beta(t,\delta)}{2N_P(t)}}.\]

Assume that $(P_{t+1} = (1,1))$. By definition of the sampling strategy, one has $\LCB_{(1,1)}(t) \leq \LCB_{(1,2)}(t)$ and $N_{(1,1)}(t) \leq N_{Y_t}(t)$. If $(\tau_\delta >t)$, one has 
\begin{eqnarray*}
 \LCB_{(1,1)}(t) & \leq & \UCB_{Y_t}(t) \\
  \UCB_{(1,1)}(t) - 2 \sqrt{\frac{\beta(t,\delta)}{2N_{(1,1)}(t)}} & \leq &  \LCB_{Y_t}(t) + 2 \sqrt{\frac{\beta(t,\delta)}{2N_{Y_t}(t)}}.
\end{eqnarray*}
On $\cE$, $\mu_{1,1} \leq \UCB_{(1,1)}(t)$ and $\LCB_{Y_t}(t) = \min (\LCB_{(2,1)}(t),\LCB_{(2,2)}(t)) \leq \min(\mu_{2,1},\mu_{2,2}) =  \mu_{2,1}$. Thus 
\begin{eqnarray*}
\mu_{1,1} - \mu_{2,1}  & \leq &  2 \sqrt{\frac{\beta(t,\delta)}{2N_{Y_t}(t)}} + 2 \sqrt{\frac{\beta(t,\delta)}{2N_{(1,1)}(t)}} \leq 4 \sqrt{\frac{\beta(t,\delta)}{2N_{(1,1)}(t)}}, \\
\end{eqnarray*}
using that $N_{(1,1)}(t) \leq N_{Y_t}(t)$. This proves that 
\[(P_{t+1} = (1,1))\cap (\tau_\delta >t) \subseteq \left(N_{(1,1)}(t) \leq \frac{8\beta(t,\delta)}{(\mu_{1,1} - \mu_{2,1})^2}\right).\]
A very similar reasoning shows that
\[(P_{t+1} = (1,2))\cap (\tau_\delta >t) \subseteq \left(N_{(1,2)}(t) \leq \frac{8\beta(t,\delta)}{(\mu_{1,2} - \mu_{2,1})^2}\right).\]

Assume that $(P_{t+1} = (2,1))$. If $(\tau_\delta >t)$, one has 
\begin{eqnarray*}
 \LCB_{X_t}(t) & \leq & \UCB_{(2,1)}(t) \\
  \UCB_{X_t}(t) - 2 \sqrt{\frac{\beta(t,\delta)}{2N_{X_t}(t)}} & \leq &  \LCB_{(2,1)}(t) + 2 \sqrt{\frac{\beta(t,\delta)}{2N_{(2,1)}(t)}}.
\end{eqnarray*}
On $\cE$, $\mu_{1,1}\leq\mu_{X_t} \leq \UCB_{X_t}(t)$ and $\LCB_{(2,1)}(t) \leq \mu_{2,1}$. Thus 
\begin{eqnarray*}
\mu_{1,1} - \mu_{2,1}  & \leq &  2 \sqrt{\frac{\beta(t,\delta)}{2N_{X_t}(t)}} + 2 \sqrt{\frac{\beta(t,\delta)}{2N_{(2,1)}(t)}} \leq 4 \sqrt{\frac{\beta(t,\delta)}{2N_{(2,1)}(t)}}, \\
\end{eqnarray*}
using that $N_{(2,1)}(t) \leq N_{X_t}(t)$. This proves that 
\[(P_{t+1} = (2,1))\cap (\tau_\delta >t) \subseteq \left(N_{(2,1)}(t) \leq \frac{8\beta(t,\delta)}{(\mu_{1,1} - \mu_{2,1})^2}\right).\]

Assume that $(P_{t+1} = (2,2))$. First, using the fact that $L_{(2,2)}(t) \leq L_{(2,1)}(t)$ yields, on $\cE$, 
\begin{eqnarray*}
 U_{(2,2)}(t) - 2 \sqrt{\frac{\beta(t,\delta)}{2N_{(2,2)}(t)}} &\leq& \mu_{2,1} \\
\mu_{2,2} - \mu_{2,1} & \leq & 2 \sqrt{\frac{\beta(t,\delta)}{2N_{(2,2)}(t)}},
\end{eqnarray*}
which leads to $N_{(2,2)}(t) \leq {2\beta(t,\delta)}/{(\mu_{2,2} - \mu_{2,1})^2}$. Then, if $(\tau_\delta >t)$, on $\cE$ (using also that $L_{(2,2)}(t) \leq L_{(2,1)}(t)$), 
\begin{eqnarray*}
 \LCB_{X_t}(t) & \leq & \UCB_{(2,2)}(t) \\
  \UCB_{X_t}(t) - 2 \sqrt{\frac{\beta(t,\delta)}{2N_{X_t}(t)}} & \leq &  \LCB_{(2,2)}(t) + 2 \sqrt{\frac{\beta(t,\delta)}{2N_{(2,2)}(t)}} \\
  \UCB_{X_t}(t) - 2 \sqrt{\frac{\beta(t,\delta)}{2N_{X_t}(t)}} & \leq &  \LCB_{(2,1)}(t) + 2 \sqrt{\frac{\beta(t,\delta)}{2N_{(2,2)}(t)}} \\  
  \mu_{1,1} - 2 \sqrt{\frac{\beta(t,\delta)}{2N_{X_t}(t)}}& \leq &  \mu_{2,1} + 2 \sqrt{\frac{\beta(t,\delta)}{2N_{(2,2)}(t)}}\\
  \mu_{1,1} - \mu_{2,1} & \leq & 4\sqrt{\frac{\beta(t,\delta)}{2N_{(2,2)}(t)}}.
\end{eqnarray*}
Thus, if $\mu_{2,2} < \mu_{1,1}$, one also has $N_{(2,2)}(t) \leq 8\beta(t,\delta)/(\mu_{1,1} - \mu_{2,1})^2$. Combining the two bounds yield 
\[(P_{t+1} = (2,2))\cap (\tau_\delta >t) \subseteq \left(N_{(2,2)}(t) \leq \frac{8\beta(t,\delta)}{\max\left(4(\mu_{2,2} - \mu_{2,1})^2,(\mu_{1,1} -\mu_{2,1})^2\right)}\right).\]

\section{Analysis of the Maximin-Racing algorithm}

\subsection{Proof of Lemma~\ref{lem:PACRacing}.}\label{proof:PACRacing}  First note that for every $P\in \cP$, introducing an i.i.d. sequence of successive observations from arm $P$, the sequence of associated empirical means $(\hat{\mu}_{P}(r))_{r \in \N}$ is defined independently of the arm being active.

We introduce the event $\cE = \cE_1 \cap \cE_2$ with
\begin{eqnarray*}
 \cE_1 & = &  \bigcap_{i =1}^K\bigcap_{\substack{(i,j)\in \Arms_i: \\ \mu_{i,j} = \mu_{i,1}}} \bigcap_{\substack{(i,j')\in \Arms_i: \\ \mu_{i,j'} > \mu_{i,1}}} \left(\forall r \in \N, f(\hat{\mu}_{i,j}(r),\hat{\mu}_{i,j'}(r)) \leq \beta(r,\delta)\right) \\
  \cE_2 &=& \bigcap_{\substack{i\in \{1,\dots,K\}: \\ \mu_{i,1} < \mu_{1,1}}} \bigcap_{\substack{(i,j)\in A_i: \\ \mu_{i,j} = \mu_{i,1}}}\bigcap_{\substack{i'\in \{1,\dots,K\}: \\ \mu_{i',1} = \mu_{1,1}}}\bigcap_{(i',j')\in A_{i'}}\left(\forall r \in \N, rf(\hat{\mu}_{i,j}(r) , \hat{\mu}_{i',j'}(r)) \leq \beta(r,\delta)\right)
\end{eqnarray*}
and the event 
\[\cF = \bigcap_{P \in \cP}\left(|\hat{\mu}_{P}(r_0)  - \mu_P | \leq \frac{\epsilon}{2}\right).\]
From \eqref{inequ:BBlock} and a union bound, $\bP(\cE^c) \leq \delta / 2$. From Hoeffding inequality and a union bound, using also the definition of $r_0$, one has $\bP(\cF^c) \leq \delta/2$. Finally, $\bP_{\bm\mu}\left(\cE \cap \cF\right) \geq 1 - \delta$.

We now show that on $\cE \cap \cF$, the algorithm outputs an $\epsilon$-optimal arm. On the event $\cE$, the following two statements are true for any round $r \leq r_0$: 
\begin{enumerate}
 \item For all $i$, if $\cR_i \neq \O$, then there exists $(i,j) \in \cR_i$ such that $\mu_{i,j}=\mu_{i,1}$ 
 \item If there exists $i$ such that $\cR_i \neq \O$, then there exists $i': \mu_{i',1}=\mu_{1,1}$ such that $\cR_{i'} \neq \O$.
\end{enumerate}
Indeed, if 1. is not true, there is a non empty set $\cR_i$ in which all the arms in the set $\{ (i,j) \in \Arms_i : \mu_{i,j}=\mu_{i,1}\}$ have been discarded. Hence, in a previous round at least one of these arms must have appeared strictly larger  than one of the arms in the set $\{ (i,j') \in \Arms_i : \mu_{i,j'}>\mu_{i,1}\}$ (in the sense of our elimination rule), which is not possible from the definition of $\cE_1$. Now if 2. is not true, there exists $i' : \mu_{i',1}=\mu_{1,1}$, such that $\cR_{i'}$ has been discarded at a previous round by some non-empty set $\cR_i$, with $\mu_{i,1}<\mu_{1,1}$. Hence, there exists $(i',j') \in A_{i'}$ that appears significantly smaller than all arms in $\cR_i$ (in the sense of our elimination rule). As $\cR_i$ contains by 1. some arm $\mu_{i,j}$ with $\mu_{i,j}=\mu_{i,1}$, there exists $r$ such that $rd(\mu_{(i,j)}(r) , \mu_{(i',j')}(r)) > \beta(r,\delta)$, which contradicts the definition of $\cE_2$.       

From the statements 1. and 2., on $\cE\cap \cF$ if the algorithm terminates before $r_0$, using that the last set in the race $\cR_i$ must satisfy $\mu_{i,1}=\mu_{1,1}$, the action $\ihat$ is in particular $\epsilon$-optimal. If the algorithm has not stopped at $r_0$, the arm $\ihat$ recommended is the empirical maximin action. Letting $\cR_i$  some set still in the race with $\mu_{i,1}=\mu_{1,1}$, one has, 
\[\min_{P \in \cR_{\ihat}} \hat{\mu}_{P}(r_0) \geq \min_{P \in \cR_i} \hat{\mu}_{P}(r_0).\]
As $\cF$ holds and because there exists $(\ihat,\jhat) \in \cR_{\ihat}$ with $\mu_{\ihat,\jhat}=\mu_{\ihat,1}$, and $(i,j) \in \cR_i$ with $\mu_{i,j}=\mu_{1,1}$, one has 
\begin{eqnarray*}
\min_{P \in \cR_i} \hat{\mu}_{P}(r_0) & \geq  &\min_{P \in \cR_i} ({\mu}_{P} - \epsilon /2)  =  \mu_{i,j} - \epsilon/2 = \mu_{1,1} - \epsilon/2.\\
 \min_{P \in \cR_{\ihat}} \hat{\mu}_{P}(r_0) & \leq & \min_{P \in \cR_{\ihat}} ({\mu}_{P} + \epsilon /2)  = \mu_{\ihat,\jhat} + \epsilon/2=\mu_{\ihat,1} + \epsilon/2.
 \end{eqnarray*}
and thus $\ihat$ is $\epsilon$-optimal, since 
\[\mu_{\ihat,1} + \frac{\epsilon}{2} \geq \mu_{1,1} - \frac{\epsilon}{2} \ \ \Leftrightarrow \ \ \mu_{1,1} - \mu_{\ihat,1} \leq \epsilon.\]
\qed 

\subsection{Proof of Theorem~\ref{thm:PACRacing}}

Recall $\mu_{1,1} > \mu_{2,1}$. We present the proof assuming additionally that for all $i\in \{1,K\}$, $\mu_{i,1} < \mu_{i,2}$ (an assumption that can be relaxed, at the cost of more complex notations). 

Let $\alpha>0$. The function $f$ defined in \eqref{def:EliminationRule} is uniformly continuous on $[0,1]^2$, thus there exists $\eta^\alpha$ such that 
\[ ||(x,y)- (x',y')||_{\infty}\leq \eta^\alpha \ \ \Rightarrow \ \ |f(x,y) - f(x',y')| \leq \alpha. \]
We introduce the event 
\[\cG_{\alpha,r} = \bigcap_{P \in \cP} (|\hat{\mu}_{P}(r) - \mu_P|\leq \eta^\alpha)\]
and let $\cE$ be the event defined in the proof of Lemma~\ref{lem:PACRacing}, which rewrites in a simpler way with our assumptions on the arms : 
\[\cE = \bigcap_{i=2}^K\bigcap_{j=1}^{K_1}\left(\forall r \in \N, rf(\hat{\mu}_{i,1}(r) , \hat{\mu}_{1,j}(r)) \leq \beta(r,\delta)\right) \bigcap_{i=1}^K\bigcap_{j=2}^{K_i} \left(\forall r \in \N, f(\hat{\mu}_{i,1}(r),\hat{\mu}_{i,j}(r)) \leq \beta(r,\delta)\right) \]
Recall that on this event, arm (1,1) is never eliminated before the algorithm stops and whenever an arm $(i,j) \in \cR$, we know that the corresponding minimal arm $(i,1)\in \cR$.

Let $(i,j) \neq (1,1)$ and recall that $\tau_\delta(i,j)$ is the number of rounds during which arm $(i,j)$ is drawn. One has 
\[\bE_{\bm\mu}[\tau_\delta(i,j)] = \bE_{\bm\mu}[\tau_\delta(i,j)\ind_{\cE}] + \bE_{\bm\mu}[\tau_\delta(i,j)\ind_{\cE^c}] \leq  \bE_{\bm\mu}[\tau_\delta(i,j)\ind_{\cE}] + \frac{r_0 \delta}{2}.\]
On the event $\cE$, if arm $(i,j)$ is still in the race at the end of round $r$,
\begin{itemize}
 \item it cannot be significantly larger than $(i,1)$: $rf(\hat{\mu}_{i,j}(r),\hat{\mu}_{i,1}(r)) \leq \beta(r,\delta)$
 \item arm $(i,1)$ cannot be significantly smaller than $(1,1)$ (otherwise all arms in $\cR_i$, including $(i,j)$, are eliminated): $rf(\hat{\mu}_{i,1}(r),\hat{\mu}_{1,1}(r)) \leq \beta(r,\delta)$ 
\end{itemize}
Finally, one can write 
\begin{eqnarray*}
 \bE_{\bm\mu}[\tau_\delta(i,j)\ind_{\cE}] & \leq & \bE_{\bm\mu}\left[\ind_{\cE} \sum_{r=1}^{r_0} \ind_{((i,j) \in \cR \ \text{at round r})} \right] \\ 
 &\leq & 
  \bE_{\bm\mu}\left[\sum_{r=1}^{r_0} \ind_{(r \max \left[f(\hat{\mu}_{i,j}(r),\hat{\mu}_{i,1}(r)), f(\hat{\mu}_{i,1}(r),\hat{\mu}_{1,1}(r))\right] \leq \beta(r,\delta) )} \right] \\
  & \leq & \bE_{\bm\mu}\left[\sum_{r=1}^{r_0} \ind_{(r \max \left[f(\hat{\mu}_{i,j}(r),\hat{\mu}_{i,1}(r)), f(\hat{\mu}_{i,1}(r),\hat{\mu}_{1,1}(r))\right] \leq \beta(r,\delta) )} \ind_{\cG_{\alpha,r}}\right] + \sum_{r=1}^{r_0} \bP_{\bm\mu}(\cG_{\alpha,r}^c) \\
  & \leq & \sum_{r=1}^{r_0} \ind_{(r \left(\max \left[f({\mu}_{i,j},{\mu}_{i,1}), f({\mu}_{i,1},{\mu}_{1,1})\right]- \alpha\right)\leq \log\left({4C_Kr}/{\delta}\right)}   + \sum_{r=1}^{\infty} \bP_{\bm\mu}(\cG_{\alpha,r}^c) \\
  & \leq & T_{(i,j)}(\delta, \alpha) + \sum_{r=1}^\infty  2\overline{K}\exp(-2(\eta^\alpha)^2 r),
\end{eqnarray*}
using Hoeffding inequality and introducing 
\[T_{(i,j)}(\delta,\alpha) : = \inf \left\{ r \in \N : r \left(\max \left[f({\mu}_{i,j},{\mu}_{i,1}), f({\mu}_{i,1},{\mu}_{1,1})\right]- \alpha\right)> \log\left(\frac{4C_Kr}{\delta}\right)\right\}\]
Some algebra (Lemma~\ref{lem:technical}) shows that $ T_{(i,j)}(\delta,\alpha) = \frac{1}{\max \left[f({\mu}_{i,j},{\mu}_{i,1}), f({\mu}_{i,1},{\mu}_{1,1})\right]- \alpha} \log\left(\frac{4C_K}{\delta}\right) + o_{\delta \rightarrow 0} \left(\log \frac{1}{\delta}\right)$ and finally, for all $\alpha >0$, 
\[\bE_{\bm\mu}\left[\tau_\delta(i,j)\right] \leq \frac{1}{\max[f({\mu}_{i,j},{\mu}_{i,1}),f({\mu}_{i,1},{\mu}_{1,1})]- \alpha} \log \left(\frac{4C_K}{\delta}\right) + o\left(\log \frac{1}{\delta}\right).\]
As this holds for all $\alpha$, and keeping in mind the trivial bound $\bE_{\bm\mu}\left[\tau_\delta(i,j)\right] \leq r_0 = \frac{2}{\epsilon^2}\log \left(\frac{4\overline{K}}{\delta}\right)$, one obtains 
\[\limsup_{\delta \rightarrow 0} \frac{\bE_{\bm\mu}[\tau_\delta(i,j)]}{\log(1/\delta)} \leq \frac{1}{\max\left[\epsilon^2/2,I_*({\mu}_{i,j},{\mu}_{i,1}),I_*({\mu}_{i,1},{\mu}_{1,1})\right] }. \]

To upper bound the number of draws of the arm $(1,1)$, one can proceed similarly and write that, for all $\alpha >0$, 
\begin{eqnarray*}
 \tau_\delta(1,1)\ind_{\cE} &=& \sup_{(i,j) \in \cP \backslash \{(1,1)\}} \tau_\delta(i,j) \ind_{\cE} \\ & \leq &  \sup_{(i,j) \in \cP \backslash \{(1,1)\}} \sum_{r=1}^{r_0} \ind_{(r \max \left[f(\hat{\mu}_{i,j}(r),\hat{\mu}_{i,1}(r)), f(\hat{\mu}_{i,1}(r),\hat{\mu}_{1,1}(r))\right] \leq \beta(r,\delta))} \\
 & \leq &  \sup_{(i,j) \in \cP \backslash \{(1,1)\}} \sum_{r=1}^{r_0} \ind_{(r (f({\mu}_{i,j},{\mu}_{i,1})\wedge f({\mu}_{i,1},{\mu}_{1,1}) - \alpha ) \leq \beta(r,\delta))} + \sum_{r=1}^{\infty} \ind_{\cG_{\alpha,r}^c} \\
 & \leq& \sup_{(i,j) \in \cP \backslash \{(1,1)\}} T_{(i,j)}(\delta,\alpha) +  \sum_{r=1}^{\infty}\ind_{\cG_{\alpha,r}^c}.
\end{eqnarray*}
Taking the expectation and using the more explicit expression of the $T_{(i,j)}$ yields 
\[\limsup_{\delta \rightarrow 0} \frac{\bE_{\bm\mu}[\tau_\delta(1,1)]}{\log(1/\delta)} \leq \frac{1}{\max\left[\epsilon^2/2,I_*({\mu}_{(2,1)},{\mu}_{1,1})\right]}. \]

\end{document}